\documentclass[12pt,reqno]{amsart}
\usepackage{amsmath, amsthm, amscd, amsfonts, amssymb, graphicx, color}
\usepackage[bookmarksnumbered, colorlinks, plainpages]{hyperref}

\theoremstyle{definition}

\theoremstyle{remark}

\numberwithin{equation}{section}
\begin{document}
\begin{center}
{\textbf{A Kenmotsu metric as a conformal $\eta$-Einstein soliton }}
\end{center}
\vskip 0.3cm
\begin{center}By\end{center}\vskip 0.3cm
\begin{center}
{Soumendu Roy \footnote{The first author is the corresponding author, supported by Swami Vivekananda Merit Cum Means Scholarship, Government of West Bengal, India.}, Santu Dey $^2$ and Arindam~~Bhattacharyya $^3$}
\end{center}
\vskip 0.3cm
\address[Soumendu Roy]{Department of Mathematics,Jadavpur University, Kolkata-700032, India}
\email{soumendu1103mtma@gmail.com}

\address[Santu Dey]{Department of Mathematics, Bidhan Chandra College, Asansol - 4, West Bengal-713303 , India}
\email{santu.mathju@gmail.com}

\address[Arindam Bhattacharyya]{Department of Mathematics,Jadavpur University, Kolkata-700032, India}
\email{bhattachar1968@yahoo.co.in}
\vskip 0.5cm
\begin{center}
\textbf{Abstract}\end{center}
The object of the present paper is to study some properties of Kenmotsu manifold whose metric is conformal $\eta$-Einstein soliton. We have studied some certain properties of Kenmotsu manifold admitting conformal $\eta$-Einstein soliton. We have also constructed a 3-dimensional Kenmotsu manifold satisfying conformal $\eta$-Einstein soliton.
\vspace {0.3cm}
\section{\textbf{Introduction}}
The notion of Einstein soliton was introduced by G. Catino and L. Mazzieri \cite{CATINO} in 2016, which generates self-similar solutions to Einstein flow,
\begin{equation}\label{1.1}
  \frac{\partial g}{\partial t} = -2(S-\frac{r}{2}g),
\end{equation}
where $S$ is Ricci tensor, $g$ is Riemannian metric and $r$ is the scalar curvature.\\
The equation of the $\eta$-Einstein soliton \cite{blaga} is given by,
\begin{equation}\label{1.2}
  \pounds_\xi g+2S+(2\lambda-r)g+2\mu \eta \otimes \eta=0,
\end{equation}
where $\pounds_\xi$ is the Lie derivative along the vector field $\xi$, $S$ is the Ricci tensor, $r$ is the scalar curvature of the Riemannian metric $g$, and $\lambda$ and $\mu$ are real constants. For $\mu=0$, the data $(g, \xi, \lambda)$ is called Einstein soliton.\\
In 2018, Mohd Danish Siddiqi \cite{mohd} introduced the notion of conformal $\eta$-Ricci soliton \cite{roy} as:
\begin{equation}\label{1.3}
  \pounds_\xi g + 2S +[2\lambda - (p + \frac{2}{n})]g+2 \mu \eta \otimes \eta=0,
\end{equation}
where  $\pounds_\xi$ is the Lie derivative along the vector field $\xi$ , $S$ is the Ricci tensor, $\lambda$, $\mu$ are constants, $p$ is a scalar non-dynamical field(time dependent scalar field)and $n$ is the dimension of manifold. For $\mu = 0$, conformal $\eta$-Ricci soliton becomes conformal Ricci soliton \cite{soumendu}.  \\
In \cite{roy1}, Roy, Dey and Bhattacharyya have defined conformal Einstein soliton, which can be written as:
\begin{equation}\label{1.3new}
  \pounds_V g+2S+[2\lambda-r+(p+\frac{2}{n})]g=0,
\end{equation}
where  $\pounds_V$ is the Lie derivative along the vector field $V$ , $S$ is the Ricci tensor, $r$ is the scalar curvature of the Riemannian metric $g$, $\lambda$ is real constant, $p$ is a scalar non-dynamical field(time dependent scalar field)and $n$ is the dimension of manifold.\\
So we introduce the notion of conformal $\eta$-Einstein soliton as:\\
\textbf{Definition 1.1:} A Riemannian manifold $(M,g)$ of dimension $n$ is said to admit conformal $\eta$-Einstein soliton if
\begin{equation}\label{1.4}
   \pounds_\xi g+2S+[2\lambda-r+(p+\frac{2}{n})]g+2\mu \eta \otimes \eta=0,
\end{equation}
where  $\pounds_\xi$ is the Lie derivative along the vector field $\xi$ ,$\lambda$, $\mu$ are real contants and  $S$, $r$, $p$, $n$ are same as defined in \eqref{1.3new}.\\\\
In the present paper we study conformal $\eta$-Einstein soliton on Kenmotsu manifold. The paper is organized as follows:\\
After introduction, section 2 is devoted for preliminaries on (2$n$+1) dimensional Kenmotsu manifold. In section 3, we have studied conformal $\eta$-Einstein soliton on Kenmotsu manifold. Here we proved if a (2$n$+1) dimensional Kenmotsu manifold admits conformal $\eta$-Einstein soliton then the manifold becomes $\eta$-Einstein. We have also characterized the nature of the manifold if the manifold is Ricci symmetric and the Ricci tensor is $\eta$-recurrent. Also we have discussed about the condition when the manifold has cyclic Ricci tensor. Then we have obtained the conditons in a (2$n$+1) dimensional Kenmotsu manifold admitting Conformal $\eta$-Einstein soliton when a vector field $V$ is pointwise co-linear with $\xi$ and a (0,2) tensor field $h$ is parallel with respect to the Levi-Civita connection associated to $g$. We have also examined the nature of a Ricci-recurrent Kenmotsu manifold admitting conformal $\eta$-Einstein soliton.\\
In last section we have given an example of a 3-dimensional Kenmotsu manifold satisfying conformal $\eta$-Einstein soliton.\\
\vspace {0.3cm}
\section{\textbf{Preliminaries}}
Let $M$ be a (2$n$+1) dimensional connected almost contact metric manifold with an almost contact metric structure $(\phi, \xi, \eta, g)$ where $\phi$ is a $(1,1)$ tensor field, $\xi$ is a vector field, $\eta$ is a 1-form  and $g$ is the compatible Riemannian metric such that\\
\begin{equation}\label{2.1}
\phi^2(X) = -X + \eta(X)\xi, \eta(\xi) = 1, \eta \circ \phi = 0, \phi \xi = 0,
\end{equation}\\
\begin{equation}\label{2.2}
g(\phi X,\phi Y) = g(X,Y) - \eta(X)\eta(Y),
\end{equation}\\
\begin{equation}\label{2.3}
g(X,\phi Y) = -g(\phi X,Y),
\end{equation}\\
\begin{equation}\label{2.4}
g(X,\xi) = \eta(X),
\end{equation}\\
for all vector fields $X, Y \in \chi(M).$\\
An almost contact metric manifold is said to be a Kenmotsu manifold \cite{kenmotsu} if
\begin{equation}\label{2.5}
  (\nabla_X \phi)Y= -g(X,\phi Y)\xi-\eta(Y)\phi X,
\end{equation}
\begin{equation}\label{2.6}
  \nabla_X \xi=X-\eta(X)\xi,
\end{equation}
where $\nabla$ denotes the Riemannian connection of $g$.\\
In a Kenmotsu manifold the following relations hold \cite{bagewadi}:
\begin{equation}\label{2.7}
  \eta(R(X,Y)Z)= g(X,Z)\eta(Y)-g(Y,Z)\eta(X),
\end{equation}
\begin{equation}\label{2.8}
  R(X,Y)\xi=\eta(X)Y-\eta(Y)X,
\end{equation}
\begin{equation}\label{2.9}
  R(X,\xi)Y=g(X,Y)\xi-\eta(Y)X,
\end{equation}
where $R$ is the Riemannian curvature tensor.\\
\begin{equation}\label{2.10}
  S(X,\xi)=-2n\eta(X),
\end{equation}
\begin{equation}\label{2.11}
  S(\phi X,\phi Y)= S(X,Y)+2n\eta(X)\eta(Y),
\end{equation}
\begin{equation}\label{2.12}
  (\nabla_X \eta)Y= g(X,Y)-\eta(X)\eta(Y),
\end{equation}
for all vector fields $X, Y, Z \in \chi(M).$\\
Now we know,
\begin{equation}\label{2.13}
  (\pounds_\xi g)(X,Y)=g(\nabla_X \xi, Y)+g(X, \nabla_Y \xi),
\end{equation}
for all vector fields $X, Y, \in \chi(M).$\\
Then using \eqref{2.6} and \eqref{2.13}, we get,
\begin{equation}\label{2.14}
  (\pounds_\xi g)(X,Y)=2[g(X,Y)-\eta(X)\eta(Y)].
\end{equation}
\vspace {0.3cm}
\section{\textbf{Conformal $\eta$-Einstein soliton on Kenmotsu manifold}}
Let $M$ be a (2$n$+1) dimensional Kenmotsu manifold. Consider the conformal $\eta$-Einstein soliton \eqref{1.4} on $M$ as:
\begin{equation}\label{3.1}
  (\pounds_\xi g)(X,Y)+2S(X,Y)+[2\lambda-r+(p+\frac{2}{2n+1})]g(X,Y)+2\mu \eta(X)\eta(Y)=0,
\end{equation}
for all vector fields $X, Y, \in \chi(M).$\\
Then using \eqref{2.14},the above equation becomes,
\begin{equation}\label{3.2}
  S(X,Y)=-[\lambda-\frac{r}{2}+\frac{(p+\frac{2}{2n+1})}{2}+1]g(X,Y)-(\mu-1)\eta(X)\eta(Y).
\end{equation}
Taking $Y=\xi$ in the above equation and using \eqref{2.10}, we get,
\begin{equation}\label{3.3}
  r=(p+\frac{2}{2n+1})-4n+2\lambda+2\mu,
\end{equation}
since $\eta(X) \neq 0$, for all $X \in \chi(M)$.\\
Also from \eqref{3.2}, it follows that the manifold is $\eta$- Einstein.\\
This leads to the following:\\\\
\textbf{Theorem 3.1.} {\em If the metric of a (2$n$+1) dimensional Kenmotsu manifold is a conformal $\eta$-Einstein soliton then the manifold becomes $\eta$- Einstein and the scalar curvature is $(p+\frac{2}{2n+1})-4n+2\lambda+2\mu$.} \\\\
We know,
\begin{equation}\label{3.4}
  (\nabla_X S)(Y,Z)= X(S(Y,Z))-S(\nabla_X Y,Z)-S(Y,\nabla_X Z),
\end{equation}
for all vector fields $X, Y, Z$ on $M$ and $\nabla$ is the Levi-Civita  connection associated with $g$.\\
Now replacing the expression of S from \eqref{3.2}, we obtain,
\begin{equation}\label{3.5}
  (\nabla_X S)(Y,Z) = -(\mu-1)[\eta(Z)(\nabla_X \eta)Y+\eta(Y)(\nabla_X \eta)Z].
\end{equation}
for all vector fields $X, Y, Z$ on $M$.\\\\
Let the manifold $M$ be Ricci symmetric i.e $\nabla S = 0.$\\
Then from \eqref{3.5}, we get,
\begin{equation}\label{3.6}
  -(\mu-1)[\eta(Z)(\nabla_X \eta)Y+\eta(Y)(\nabla_X \eta)Z]=0,
\end{equation}
for all vector fields $X, Y, Z \in \chi(M).$\\
Taking $Z = \xi$ in the above equation and using \eqref{2.12}, \eqref{2.1}, we obtain,
\begin{equation}\label{3.7}
  \mu=1.
\end{equation}
Then from \eqref{3.3}, we get,
\begin{equation}\label{3.8}
  r=(p+\frac{2}{2n+1})-4n+2\lambda+2.
\end{equation}
So we can state the following theorem:\\\\
\textbf{Theorem 3.2.} {\em If the metric of a (2$n$+1) dimensional Ricci symmteric Kenmotsu manifold is a conformal $\eta$-Einstein soliton then $\mu$ = 1 and the scalar curvature is $(p+\frac{2}{2n+1})-4n+2\lambda+2$.}\\\\
Now if the Ricci tensor $S$ is $\eta$-recurrent, then we have,
\begin{equation}\label{3.9}
  \nabla S= \eta \otimes S,
\end{equation}
which implies that,
\begin{equation}\label{3.10}
  (\nabla_X S)(Y,Z)=\eta(X)S(Y,Z),
\end{equation}
for all vector fields $X, Y, Z$ on $M$.\\
Using \eqref{3.5}, the above equation reduces to,
\begin{equation}\label{3.11}
  -(\mu-1)[\eta(Z)(\nabla_X \eta)Y+\eta(Y)(\nabla_X \eta)Z]=\eta(X)S(Y,Z).
\end{equation}
Taking $Y = \xi, Z = \xi$ in the above equation and using \eqref{2.12},\eqref{3.2},  we get,
\begin{equation}\label{3.12}
  [\lambda+\mu-\frac{r}{2}+\frac{p+\frac{2}{2n+1}}{2}]\eta(X)=0,
\end{equation}
which implies that,
\begin{equation}\label{3.13}
  r=2\lambda+2\mu+(p+\frac{2}{2n+1}).\\
\end{equation}
Then we can state the following:\\\\
\textbf{Theorem 3.3} {\em If the metric of a (2$n$+1) dimensional Kenmotsu manifold is a conformal $\eta$-Einstein soliton and the Ricci tensor $S$ is $\eta$- Recurrent, then the scalar curvature is $2\lambda + 2\mu+ (p+\frac{2}{2n+1})$.}\\\\
Similarly from \eqref{3.5}, we get,
\begin{equation}\label{3.14}
  (\nabla_Y S)(Z,X) = -(\mu-1)[\eta(X)(\nabla_Y \eta)Z+\eta(Z)(\nabla_X \eta)Y],
\end{equation}
and
\begin{equation}\label{3.15}
  (\nabla_Z S)(X,Y) = -(\mu-1)[\eta(Y)(\nabla_Z \eta)X+\eta(X)(\nabla_Z \eta)Y].
\end{equation}
for all vector fields $X, Y, Z$ on $M$.\\
Then adding \eqref{3.5},\eqref{3.14}, \eqref{3.15} and using \eqref{2.12}, \eqref{2.2}, we obtain,
\begin{eqnarray}\label{3.16}
  (\nabla_X S)(Y,Z)+(\nabla_Y S)(Z,X)+(\nabla_Z S)(X,Y) &=& -2(\mu-1)[\eta(X)g(\phi Y,\phi Z) \nonumber \\
  &+& \eta(Y)g(\phi Z,\phi X) \nonumber\\
  &+& \eta(Z)g(\phi X,\phi Y)].
\end{eqnarray}\\
Now if the manifold $M$ has  cyclic Ricci tensor i.e $(\nabla_X S)(Y,Z) + (\nabla_Y S)(Z,X) + (\nabla_Z S)(X,Y) = 0,$ then from \eqref{3.16}, we have,
\begin{equation}\label{3.17}
  (\mu-1)[\eta(X)g(\phi Y,\phi Z)+\eta(Y)g(\phi Z,\phi X)+\eta(Z)g(\phi X,\phi Y)]=0.
\end{equation}
Taking $X = \xi$ in the above equation and using \eqref{2.1}, we get,
\begin{equation}\label{3.18}
  \mu=1.
\end{equation}
Again if we take $\mu = 1$ in \eqref{3.16}, we obtaion $(\nabla_X S)(Y,Z) + (\nabla_Y S)(Z,X) + (\nabla_Z S)(X,Y) = 0,$ i.e the manifold $M$ has cyclic Ricci tensor.\\
Hence we can state the following:\\\\
\textbf{Theorem 3.4} {\em Let the metric of a (2$n$+1) dimensional Kenmotsu manifold $M$ is a conformal $\eta$-Einstein soliton. Then $M$ has cyclic Ricci tensor iff $\mu = 1.$}\\\\
Now if $\mu = 1$, then from \eqref{3.3} we obtain,
\begin{equation}\label{3.19}
r=(p+\frac{2}{2n+1})-4n+2\lambda+2.
\end{equation}
Then we have,\\
\textbf{Corollary 3.5.} {\em If a (2$n$+1) dimensional Kenmotsu manifold $M$ has a cyclic Ricci tensor and the metric is a conformal $\eta$-Einstein soliton then the scalar curvature is $(p+\frac{2}{2n+1})-4n+2\lambda+2$.}\\\\
Let a conformal $\eta$-Einstein soliton is defined on a (2$n$+1) dimensional Kenmotsu manifold $M$ as,
\begin{equation}\label{3.20}
 \pounds_V g+2S+[2\lambda-r+(p+\frac{2}{2n+1})]g+2\mu \eta \otimes \eta=0,
\end{equation}
where  $\pounds_V$ is the Lie derivative along the vector field $V$ , $S$ is the Ricci tensor, $r$ is the scalar curvature of the Riemannian metric $g$, $\lambda$, $\mu$ are real contants, $p$ is a scalar non-dynamical field(time dependent scalar field).\\
Let $V$ be pointwise co-linear with $\xi$, i.e $V = b\xi$, where $b$ is a function on $M$.\\
Then \eqref{3.20} becomes,
\begin{equation}\label{3.21}
  (\pounds_{b\xi} g)(X,Y)+2S(X,Y)+[2\lambda-r+(p+\frac{2}{2n+1})]g(X,Y)+2\mu \eta(X)\eta(Y)=0,
\end{equation}
for all vector fields $X, Y$ on $M$.\\
Applying the property of Lie derivative and Levi-Civita connection we have,
\begin{multline}\label{3.22}
   bg(\nabla_X \xi,Y)+(Xb)\eta(Y)+bg(\nabla_Y \xi,X)+(Yb)\eta(X)+2S(X,Y)\\
  + [2\lambda-r+(p+\frac{2}{2n+1})]g(X,Y)+2\mu \eta(X)\eta(Y)=0.
\end{multline}
Now using \eqref{2.6}, we get,
\begin{multline}\label{3.23}
 2bg(X,Y)-2b\eta(X)\eta(Y)+(Xb)\eta(Y)+(Yb)\eta(X)+2S(X,Y) \\
  + [2\lambda-r+(p+\frac{2}{2n+1})]g(X,Y)+2\mu \eta(X)\eta(Y)=0.
\end{multline}
Taking $Y = \xi$ in the above equation and using \eqref{2.1},\eqref{2.4},\eqref{2.10}, we obtain,
\begin{equation}\label{3.24}
  (Xb)+(\xi b)\eta(X)-4n\eta(X)+[2\lambda-r+(p+\frac{2}{2n+1})]\eta(X)+2\mu \eta(X)=0.
\end{equation}
Then by putting $X = \xi$, the above equation reduces to,
\begin{equation}\label{3.25}
  \xi b=2n+\frac{r}{2}-\lambda-\mu-\frac{(p+\frac{2}{2n+1})}{2}.
\end{equation}
Using \eqref{3.25}, \eqref{3.24} becomes,
\begin{equation}\label{3.26}
  (Xb)+[\lambda+\mu+\frac{(p+\frac{2}{2n+1})}{2}-2n-\frac{r}{2}]\eta(X)=0.
\end{equation}
Applying exterior differentiation in \eqref{3.26}, we have,
\begin{equation}\label{3.27}
  [\lambda+\mu+\frac{(p+\frac{2}{2n+1})}{2}-2n-\frac{r}{2}]d\eta=0.
\end{equation}
Now we know,
\begin{equation}\label{3.28}
d\eta(X,Y)=\frac{1}{2}[(\nabla_X \eta)Y-(\nabla_Y \eta)X],
\end{equation}
for all vector fields $X, Y$ on $M$.\\
Using \eqref{2.12},the above equation becomes,
\begin{equation}\label{3.29}
  d\eta=0.
\end{equation}
Hence the 1-form $\eta$ is closed.\\
So from \eqref{3.27}, either $r = 2\lambda + 2\mu +  (p + \frac{2}{2n+1}) - 4n$ or $r \neq 2\lambda + 2\mu +  (p + \frac{2}{2n+1}) - 4n$.\\
If $r = 2\lambda + 2\mu +  (p + \frac{2}{2n+1}) - 4n$, \eqref{3.26} reduces to,
\begin{equation}\label{3.30}
  (Xb)=0.
\end{equation}
This implies that $b$ is constant.\\\\
So we can state the following theorem:\\\\
\textbf{Theorem 3.6.} {\em Let $M$ be a (2$n$+1) dimensional Kenmotsu manifold admitting a conformal $\eta$-Einstein soliton $(g,V)$, $V$ being a vector field on $M$. If $V$ is pointwise co-linear with $\xi$, a vector field on $M$, then $V$ is a constant multiple of $\xi$, provided the scalar curvature is $2\lambda + 2\mu +  (p + \frac{2}{2n+1}) - 4n$.}\\\\
Let $h$ be a symmetric tensor field of (0,2) type which we suppose to be parallel with respect to the Levi-Civita connection $\nabla$ i.e $\nabla h = 0.$\\
Applying the Ricci commutation identity, we have,
\begin{equation}\label{3.31}
  \nabla^2 h (X,Y;Z,W)-\nabla^2 h(X,Y;W,Z)=0.
\end{equation}
for all vector fields $X, Y, Z, W$ on $M$.\\
From \eqref{3.31}, we obtain the relation,
\begin{equation}\label{3.32}
  h(R(X,Y)Z,W)+h(Z,R(X,Y)W)=0.
\end{equation}
Replacing $Z = W = \xi$ in the above equation and using \eqref{2.8}, we get,
\begin{equation}\label{3.33}
  \eta(X)h(Y,\xi)-\eta(Y)h(X,\xi)=0.
\end{equation}
Replacing $X = \xi$ and using \eqref{2.1}, the above equation reduces to,
\begin{equation}\label{3.34}
  h(Y,\xi)=\eta(Y)h(\xi,\xi),
\end{equation}
for all vector fields $Y$ on $M$.\\
Differentiating the above equation covariantly with respect to $X$, we get,
\begin{equation}\label{3.35}
  \nabla_X (h(Y,\xi))=\nabla_X(\eta(Y)h(\xi,\xi)).
\end{equation}
Now expanding the above eqution by using \eqref{3.34}, \eqref{2.6},\eqref{2.12} and the property that $\nabla h = 0$, we obtain,
\begin{equation}\label{3.36}
  h(X,Y)=h(\xi,\xi)g(X,Y),
\end{equation}
for all vector fields $X, Y$ on $M$.\\
Let us take,
\begin{equation}\label{3.37}
  h = \pounds_\xi g + 2S + 2\mu \eta \otimes \eta.
\end{equation}
Then from \eqref{2.14},\eqref{3.2}, we get,
\begin{equation}\label{3.38}
  h(\xi,\xi)=-2\lambda-(p+\frac{2}{2n+1})+r.
\end{equation}
Then using \eqref{3.37}, \eqref{3.36} becomes,
\begin{equation}\label{3.39}
  (\pounds_\xi g)(X,Y)+2S(X,Y)+[2\lambda-r+(p+\frac{2}{2n+1})]g(X,Y)+2\mu \eta(X)\eta(Y)=0,
\end{equation}
which is the Conformal $\eta$-Einstein soliton.
This leads to,\\\\
\textbf{Theorem 3.7.} {\em In a (2$n$+1) dimensional Kenmotsu manifold assume that a symmetric (0,2) tensor field $h = \pounds_\xi g + 2S + 2\mu \eta \otimes \eta$ is parallel with respect to the Levi-Civita connection associated to $g$. Then $(g , \xi)$ yields a conformal $\eta$-Einstein soliton.}\\\\
\textbf{Definition 3.8} A Kenmotsu manifold is said to be Ricci-recurrent manifold if there exists a non-zero 1-form $A$ such that\\
\begin{equation}\label{3.40}
  (\nabla_W S)(Y,Z)=A(W)S(Y,Z),
\end{equation}
for any vector fields $W, Y, Z$ on $M$.\\
Replacing $Z$ by $\xi$ in the above equation and using \eqref{2.10}, we get,
\begin{equation}\label{3.41}
  (\nabla_W S)(Y,\xi)=-2n A(W)\eta(Y),
\end{equation}
which implies that,
\begin{equation}\label{3.42}
  WS(Y,\xi)-S(\nabla_W Y,\xi)-S(Y,\nabla_W \xi)=-2n A(W)\eta(Y).
\end{equation}
Using \eqref{2.10} and \eqref{2.6}, the above equation becomes,
\begin{equation}\label{3.43}
  2n(\nabla_W \eta)Y+2n\eta(W)\eta(Y)+S(Y,W)=2n A(W)\eta(Y).
\end{equation}
Again using \eqref{2.12}, the above equation reduces to,
\begin{equation}\label{3.44}
  2ng(W,Y)+S(Y,W)=2n A(W)\eta(Y).
\end{equation}
Taking $W = \xi$ in the above equation and using \eqref{3.2}, we obtain,
\begin{equation}\label{3.45}
  r=2\lambda+2\mu+(p + \frac{2}{2n+1})+4n(A(\xi)-1).
\end{equation}
So we can state,\\\\
\textbf{Theorem 3.9.} {\em If the metric of a (2$n$+1) dimensional Ricci-recurrent Kenmotsu manifold is a conformal $\eta$-Einstein soliton with the 1-form $A$, then the scalar curvature becomes $2\lambda + 2\mu+ (p+\frac{2}{2n+1})+4n(A(\xi)-1)$.}\\\\
\section{\textbf{Example of a 3-dimensional Kenmotsu manifold admitting conformal $\eta$-Einstein soliton:}}
We consider the three-dimensional manifold $M = \{(x, y, z) \in \mathbb{R}^3 , (x, y, z) \neq (0, 0, 0) \}$, where $(x, y, z)$ are standard coordinates in $\mathbb{R}^3$. The vector fields
\begin{equation}
  e_1 = z \frac{\partial}{\partial x}, \quad e_2= z \frac{\partial}{\partial y}, \quad e_3=-z \frac{\partial}{\partial z} \nonumber
\end{equation}
are linearly independent at each point of $M$. Let $g$ be the Riemannian metric defined by
\begin{equation}
  g(e_1,e_2)=g(e_2,e_3)=g(e_3,e_1)=0,\nonumber
\end{equation}
\begin{equation}
   g(e_1,e_1) = g(e_2,e_2) = g(e_3,e_3) =1.\nonumber
\end{equation}
Let $\eta$ be the 1-form defined by $\eta(Z) = g(Z,e_3)$, for any $Z \in \chi(M)$,where $\chi(M)$ is the set of all differentiable vector fields on $M$ and $\phi$ be the (1, 1)-tensor field defined by,
\begin{equation}
  \phi e_1=-e_2, \quad \phi e_2=e_1,\quad  \phi e_3=0.\nonumber
\end{equation}
Then using the linearity of $\phi$ and $g$, we have, $$\eta(e_3) = 1, \quad \phi ^2 Z = -Z + \eta(Z)e_3, \quad g(\phi Z,\phi W) = g(Z,W) - \eta(Z)\eta(W),$$ for any $Z,W \in \chi(M)$. Thus for $e_3 = \xi$, $(\phi, \xi, \eta, g)$ defines an almost contact metric structure on $M$.\\
Let $\nabla$ be the Levi-Civita connection with respect to the Riemannian metric $g$. Then we have,
  $$ [e_1,e_2] =0, \quad [e_1,e_3] = e_1, \quad [e_2,e_3]= e_2.$$
The connection $\nabla$ of the metric $g$ is given by,
\begin{eqnarray}
  2g(\nabla_X Y,Z) &=& Xg(Y,Z)+Yg(Z,X)-Zg(X,Y)\nonumber \\
                   &-& g(X, [Y,Z])-g(Y, [X, Z]) + g(Z, [X, Y ]),\nonumber
\end{eqnarray}
which is known as Koszul’s formula.\\
Using Koszul’s formula, we can easily calculate,
$$\nabla_{e_1} e_1 =-e_3, \quad \nabla_{e_1} e_2 =0 ,\quad \nabla_{e_1} e_3 =e_1,$$
$$\nabla_{e_2} e_1 =0, \quad  \nabla_{e_2} e_2 =-e_3, \quad  \nabla_{e_2} e_3 =e_2,$$
$$\nabla_{e_3} e_1 =0, \quad  \nabla_{e_3} e_2 = 0, \quad  \nabla_{e_3} e_3 =0.$$
From the above it follows that the manifold satisfies $\nabla_X \xi = X - \eta(X)\xi$, for $\xi = e_3$. Hence the manifold is a Kenmotsu Manifold.\\
Also, the Riemannian curvature tensor $R$ is given by,
$$R(X, Y )Z = \nabla_X\nabla_Y Z - \nabla_Y \nabla_X Z - \nabla_{[X,Y]} Z.$$
Hence,
$$R(e_1,e_2)e_2 =-e_1, \quad R(e_1,e_3)e_3 = -e_1, \quad R(e_2,e_1)e_1 = -e_2,$$
$$R(e_2,e_3)e_3 = -e_2, \quad R(e_3,e_1)e_1 = -e_3, \quad R(e_3,e_2)e_2 = -e_3,$$
$$R(e_1,e_2)e_3 = 0, \quad R(e_2,e_3)e_1 = 0, \quad R(e_3,e_1)e_2 = 0.$$
Then, the Ricci tensor $S$ is given by,
$$S(e_1,e_1) = -2, \quad S(e_2,e_2) = -2, \quad S(e_3,e_3)= -2.$$
From \eqref{3.2}, we have,
\begin{equation}\label{4.1}
   S(e_3,e_3)=-[\lambda+\mu-\frac{r}{2}+\frac{(p+\frac{2}{3})}{2}],
\end{equation}
which implies that,
\begin{equation}\label{4.2}
  r=2\lambda+2\mu-4+(p+\frac{2}{3}).
\end{equation}
Hence $\lambda$ and $\mu$ satisfies equation \eqref{3.3} and so $g$ defines a conformal $\eta$-Einstein soliton on the 3-dimensional Kenmotsu manifold $M$.

\end{document}